\documentclass[a4paper,11pt,reqno]{amsart}
\usepackage[centering,includeheadfoot,margin=2.5 cm]{geometry}
\usepackage{graphics,enumitem,epsfig,textcomp}
\usepackage{amsfonts,amsmath,amssymb,euscript,color,mathrsfs}
\usepackage[utf8]{inputenc}	

\usepackage[caption=false,subrefformat=parens,labelformat=parens]{subfig}
\usepackage{url}
\usepackage[symbol]{footmisc}
\usepackage{hyperref}
\hypersetup{
    colorlinks=true,
    linkcolor=blue,
    filecolor=blue,      
    urlcolor=blue,
    citecolor=blue,
    }

\usepackage[square,numbers]{natbib}
\usepackage{amsmath,amsthm,amssymb,amsfonts,lineno}

\numberwithin{equation}{section}
\newtheorem{theorem}{Theorem}%[section]

\newtheorem{lemma}{Lemma}%[section]
\newtheorem{proposition}{Proposition}%[section]
\def\div{\nabla\cdot}
\def\diver{\div}

\def\d{\,\mathrm{d}}
\def\eps{\varepsilon}
\def\R{\mathbb{R}}
\def\C{\hbox{\rlap{\kern.24em\raise.1ex\hbox
      {\vrule height1.3ex width.9pt}}C}}
\def\P{\hbox{\rlap{I}\kern.16em P}}
\def\Q{\hbox{\rlap{\kern.24em\raise.1ex\hbox
      {\vrule height1.3ex width.9pt}}Q}}
\def\M{\hbox{\rlap{I}\kern.16em\rlap{I}M}}
\def\Z{\hbox{\rlap{Z}\kern.20em Z}}
\def\({\begin{eqnarray}}
\def\){\end{eqnarray}}
\def\[{\begin{eqnarray*}}
\def\]{\end{eqnarray*}}
\def\part#1#2{\frac{\partial #1}{\partial #2}}

\def\grad{\nabla}
\def\pmb#1{\setbox0=\hbox{$#1$}
  \kern-.025em\copy0\kern-\wd0
  \kern-.05em\copy0\kern-\wd0
  \kern-.025em\raise.0433em\box0 }

\def\totk#1#2#3{{\frac{\d^#3 #1}{\d #2^#3}}}
\def\laplace{\Delta}
\def\d{\,\mathrm{d}}

\def\R{\mathbb{R}}

\def\epsilon{\varepsilon}

\def\P{\mathbb{P}}
\def\Q{\mathbb{Q}}

\newcommand{\dx}{\mathrm{d}x}

\newcommand{\D}{\mathrm{D}}
\newcommand{\z}{\mathrm{z}}

\def\cJH#1{\textcolor{blue}{\bf [#1]}}

\title{Emergence of biological transportation networks as a self-regulated process}
\date{April 2022}

\begin{document}
\pagenumbering{gobble}
\maketitle
\pagenumbering{arabic}

\centerline{
     {\large Jan Haskovec}\footnote{Mathematical and Computer Sciences
            and Engineering Division,
         King Abdullah University of Science and Technology,
         Thuwal 23955-6900, Kingdom of Saudi Arabia;
         {\it jan.haskovec@kaust.edu.sa}}\qquad
     {\large Peter Markowich}\footnote{Mathematical and Computer Sciences and Engineering Division,
         King Abdullah University of Science and Technology,
         Thuwal 23955-6900, Kingdom of Saudi Arabia;
         {\it peter.markowich@kaust.edu.sa}, and
         Faculty of Mathematics, University of Vienna,
        Oskar-Morgenstern-Platz 1, 1090 Vienna;
         {\it peter.markowich@univie.ac.at}}\qquad
         {\large Simone Portaro}\footnote{Mathematical and Computer Sciences
            and Engineering Division,
         King Abdullah University of Science and Technology,
         Thuwal 23955-6900, Kingdom of Saudi Arabia;
         {\it simone.portaro@kaust.edu.sa}}
     }
%\begin{big}
\medskip \medskip

\centerline{\textit{Dedicated to Juan Luis V{\'a}zquez at the occasion of his 75th birthday.}}
\medskip \medskip

% riga bianca da mettere

%\end{big}

\textbf{Abstract.} We study self-regulating processes modeling biological transportation networks. Firstly, we write the formal $L^2$-gradient flow for the symmetric tensor valued diffusivity $D$ of a broad class of entropy dissipations associated with a purely diffusive model. The introduction of a prescribed electric potential leads to the Fokker-Planck equation, for whose entropy dissipations we also investigate the formal $L^2$-gradient flow. We derive an integral formula for the second variation of the dissipation functional, proving convexity (in dependence of diffusivity tensor) for a quadratic entropy density modeling Joule heating. Finally, we couple in the Poisson equation for the electric potential obtaining the Poisson-Nernst-Planck system. The formal gradient flow of the associated entropy loss functional is derived, giving an evolution equation for $D$ coupled with two auxiliary elliptic PDEs.
\medskip \medskip

\medskip \medskip

\section{Introduction}
A {transportation network} is a realization of a spatial structure which permits
flow of some commodity.
Network structures and dynamics in biological contexts, in particular
organization of leaf venation networks, vascular and neural network formation,
have been widely investigated in the recent literature \cite{morel, couder, Corson, magnasco}.
One typically focuses on studying optimality in transport properties (electric, fluids, material) of the networks,
involving a complex trade-off between cost, transportation efficiency, and fault tolerance \cite{Barthelemy}.
Biological transportation networks develop without centralized control \cite{Tero}
and have been fine-tuned by many cycles of evolutionary selection pressure.
They can therefore be considered as emergent structures resulting from self-regulating processes.

An important class of self-regulating processes that we shall study in this paper,
is governed by the minimization of an entropy dissipation, coupled to 
the conservation law for a quantity $u=u(x)$,
\(   \label{eq:ell}
   %\part{u}{t}
       -\grad\cdot (D\grad u) = S \qquad \mathrm{in} \, \, \Omega, %\qquad \mbox{in } \Omega,
\)
with $D=D(x)$ the symmetric, positive definite diffusivity tensor of the transportation structure and $S=S(x)$ the distribution
of sources and sinks.
The quantity $u=u(x)$ typically represents the concentration of a chemical species, ions, nutrients or material pressure.
Equation \eqref{eq:ell} is posed on a bounded domain $\Omega\subset\R^d$, $d\geq 1$, with smooth boundary $\partial\Omega$,
subject to the Dirichlet boundary condition
\(  \label{BC_u}
   u \equiv c \qquad \mathrm{on} \, \, \partial\Omega, %\qquad \mbox{on } \partial\Omega,
\)
where $c$ is a constant, considered the equilibrium state of the system.
%and the initial condition
%  \(
%    u(t=0) = u_{0} \geq 0 \qquad\mbox{in } \Omega.
%  \)

The entropy dissipation is given by the functional
\( \label{funct0}
   E[D] = \int_\Omega \Phi''(u) \grad u \cdot D \grad u \,\dx,
\)
where the entropy (or free energy) generating function $\Phi:\R\to\R$ is convex with
\(
   \label{eq:equilibrium}
    \Phi'(c) = 0.
\)
Here the solution $u$ of \eqref{eq:ell}, \eqref{BC_u} is considered to depend on the diffusivity tensor $D$, i.e., $u = u[D]$.
%Note that $E[D]$ is obtained by multiplication of \eqref{eq:ell} by $\Phi'(u)$ and integrating by parts,
Evolutionary selection is assumed to take place through the minimization of the entropy dissipation with respect to $D$.

A generic example is the minimization of Joule heating, where $\Phi(u) = u^2/2$ and the energy dissipation is given by the Dirichlet integral
considered as a functional of the diffusivity $D$,
\(  \label{eq:Dir}
   E[D] = \int_\Omega \grad u \cdot D \grad u \,\dx,
\)
where $u=u[D]$.
Indeed, Joule's law asserts that the system power density is given by the product of the current $J = D \nabla u$ and the potential gradient $\grad u$.
Another typical choice for the entropy generator is $\Phi(u) = u(\ln(u)-1)$, which turns \eqref{funct0}
into the Fisher information $\int_\Omega \frac{\grad u \cdot D \grad u}{u} \,\dx$.

It is easy to check that the gradient flow of \eqref{eq:Dir}
subject to the constraint \eqref{eq:ell}, \eqref{BC_u} is given by
\(  \label{eq:GF_Dir}
    \part{D}{t} = \nabla u \otimes \nabla u \qquad \mathrm{in} \, \, (0, \infty) \times \Omega,  %\qquad \mbox{in } (0, \infty) \times \Omega,
\)
where $D=D(t,x)$ and $t$ is the time-like variable induced by the gradient flow.
A particular case of this type of process in the context of biological applications
(e.g., leaf venation in plants) is the network formation problem introduced in \cite{Hu-Cai} and further analyzed
in the series of papers \cite{AAFM, bookchapter, BHMR, HMP15, HMPS16, HMP19, Li, Xu, Xu2}.
Here the quantity $D=D(t,x)$ represent the tensor-valued local conductivity of the network,
which is understood as a continuous porous medium.
The flow of the material (e.q., water with nutrients in the case of leaf venation) is described
in terms of the flux $q=D\grad u$, where $u=u(t,x)$ is the fluid pressure.
To account for the metabolic cost of maintaining the biological tissue,
the functional \eqref{eq:Dir} is extended by adding the algebraic term $\int_\Omega |D|^\gamma \d x$,
where $|D|$ denotes a suitable matrix norm of $D$ and $\gamma>0$ is
the metabolic exponent derived from the biological properties of the underlying system;
see \cite{AAFM, bookchapter, Hu} for details on the modeling.
Moreover, random fluctuations in the media are accounted for by adding the
Dirichlet integral $\frac{\beta}{2}\int_{\Omega} |\grad D|^2 \,\dx$
with $\beta>0$ the diffusivity constant.
One thus arrives at the energy functional
\(   \label{E_HuCai}
   E[D] = \int_{\Omega } \frac{\beta}{2}|\grad D|^2 + \grad u \cdot D \grad u  + \frac{\alpha}{\gamma} |D|^\gamma \,\dx,
\)
where the parameter $\alpha>0$ is the metabolic coefficient.
The energy is constrained by \eqref{eq:ell}, \eqref{BC_u}, representing the
local mass conservation $-\div q = S$.
We refer to \cite{HKM, HKM2, CMS-2022} for a derivation of the system from
the discrete graph-based model of \cite{Hu-Cai}.
A crucial observation about \eqref{eq:ell}, \eqref{BC_u}, \eqref{E_HuCai} made in \cite{CMS-2022} 
is that for $\gamma\geq 1$ it is a convex functional in $D$.
Consequently, by standard theory \cite{ABS} we obtain
the existence and uniqueness of the corresponding $L^2$-gradient flow,
\( \label{eq:GF_Dir_complete}
   \part{D}{t} = \beta\laplace D + \grad u \otimes \grad u - \alpha |D|^{\gamma-2} D \qquad \mathrm{in} \, \, (0, \infty) \times \Omega, %\qquad \mbox{in } (0, \infty) \times \Omega, 
\)
subject to a homogeneous Dirichlet boundary condition for $D$ and coupled to \eqref{eq:ell}, \eqref{BC_u}.
However, let us note that well-posedness of solutions for $\gamma\in (0,1)$
is an open problem, complicated by the singularity of the term $|D|^{\gamma-2} D$
when $D\to 0$.

The first goal of the paper, carried out in Section \ref{sec:GDM},
is to derive formal $L^2$-gradient flows corresponding to general entropy dissipation functionals
of the form \eqref{funct0} with convex functions $\Phi$.
This class of functionals is a classical subject of interest in the literature studying convex and logarithmic
Sobolev inequalities and the rate of convergence to equilibrium for Fokker-Planck type equations \cite{AMT, AMTU}.

In Section \ref{sec:DD} we extend the model to the drift-diffusion setting,
i.e., we replace \eqref{eq:ell} by
\(  \label{eq:DD}
       %\part{D}{t}
          - \div( D \nabla u + \z u D \nabla \varphi ) = S \qquad \mathrm{in} \, \, \Omega. %\qquad \mbox{on } \Omega.   
\)
The Fokker-Planck equation \eqref{eq:DD} describes flow of ions with charge $\z \in \R$
under the effects of the concentration gradient $\nabla u$
and electric field $- \z \nabla \varphi$.
In certain applications, e.g. magnetofluidics \cite{Hejazian} or ferrohydrodynamics \cite{Khan, ferro},
the influence of the particle velocity and magnetic potential is relevant.
However, in the context of biological systems, where fluid flow velocities are small, these effects are mostly negligible.
In \cite{constantin2019, schmuck2009}, the Fokker-Planck equation was coupled to the Navier-Stokes equation,
modelling the velocity of charged particles in the fluid.

To account for drift induced by the electrostatic field of the charged particles, we shall consider the electric potential $\varphi=\varphi(t,x)$
to be a solution of the Poisson equation
\(  \label{eq:Poisson}
   -\laplace\varphi = \z u \qquad \mathrm{in} \, \, \Omega, %\qquad \mbox{on } \Omega,
\)
where $\z\in\R$ is the particle charge, subject to a homogeneous Dirichlet boundary condition. 
The coupled system \eqref{eq:DD}, \eqref{eq:Poisson} then constitutes the Poisson-Nernst-Planck system,
widely used in the literature to describe the flow of ions in various physical applications such as semiconductor charge carriers' transport
\cite{markowich1985, markowich1990}, ions transport in porous media \cite{gagneux2016}, in biological processes \cite{eisenberg1996, lu2011} and in electronic devices \cite{romano}.
The free energy for Poisson-Nernst-Planck system \eqref{eq:DD}, \eqref{eq:Poisson} is given by the Helmholtz free energy, see, e.g., \cite{lu2011},
\(    \label{eq:helmholtz_ddmod}
    \mathcal{H}(u, \varphi) :=  \mathcal{S}(u) + U(\varphi) =  \int_{\Omega} u (\ln u -1) \dx + \int_{\Omega} \frac{1}{2} | \nabla \varphi |^2 \dx,
\)
%with $k_{B}$ the Boltzmann constant.
%Again, we normalize the system temperature $\mathcal{T}$ and the Boltzmann constant $k_B$ to one.
Observe that $\mathcal{S}(u) = \int_{\Omega} u (\ln u -1) \dx$ is the Boltzmann entropy
and $U(\varphi) = \int_{\Omega} \frac{1}{2} | \nabla \varphi |^2 \dx$ is the internal energy due to electrostatic particle interactions.

In Section \ref{sec:PNP} we calculate the loss of the Helmholtz free energy to be of the form
\(
    \label{eq:E_PNP}
    \mathcal{E}[D] = \int_{\Omega} u \nabla \mu \cdot D \nabla \mu \, \dx.
\)
with the quasi-Fermi energy level $\mu$ given by
\(   \label{eq:mu_ddmod}
    \mu := \ln(u) + \z \varphi.
\)
Let us note that functionals of the form \eqref{eq:E_PNP} also appear in the context of charged particle flow in semiconductor devices,
see, e.g., \cite{boukili2013}.
%For the tensor-valued diffusivity $D$ the Joule heating term is given by
%\[
%    Q := \frac{J \cdot D^{-1} J}{u} = - J \cdot \nabla \mu = u \nabla \mu \cdot D \nabla \mu.
%\]

The main goal of this paper is to calculate the gradient flow of the loss functional \eqref{eq:E_PNP}
constrained by the Poisson-Nernst-Planck system \eqref{eq:DD}, \eqref{eq:Poisson}.
The calculation shall be carried out in Section \ref{sec:PNP}, leading to the evolution
equation for the diffusivity tensor $D=D(t,x)$
\(    \label{eq:gradient_flow_PNP}
    \part{D}{t} = u \nabla \mu \otimes \nabla \mu + u \frac{\nabla\mu\otimes\nabla\sigma + \nabla\sigma\otimes\nabla\mu}{2} \qquad \mathrm{in} \, \, (0, \infty) \times \Omega,
\)
coupled to the system for the auxiliary variables $\sigma=\sigma(t,x)$ and $\eta=\eta(t,x)$,
\begin{linenomath}
%\begin{subequations}
\begin{align}
%\label{eq:auxiliary_PNP1}
    - \div(D \nabla \sigma) + \z \nabla \varphi \cdot D \nabla \sigma - \z^2 \eta = \nabla \mu \cdot D \nabla \mu \qquad &\mathrm{in} \, \, (0, \infty) \times \Omega, \notag \\
    %\label{eq:auxiliary_PNP2}
    - \Delta \eta = \div(u D \nabla \sigma) \qquad &\mathrm{in} \, \, (0, \infty) \times \Omega \notag,
\end{align}
%\end{subequations}
\end{linenomath}
subject to the homogeneous Dirichlet boundary conditions for $\sigma$ and $\eta$,
\(  %\label{eq:auxiliary_BC1}
        \sigma = 0, \quad \eta = 0 \qquad \mbox{on } \partial \Omega. \notag
\)

When equipped with a metabolic term of the form $\mu\int_\Omega |D|^\gamma \, \d x$ and with diffusion, analogously to \eqref{E_HuCai},
the loss of the Helmholtz free energy \eqref{eq:E_PNP} becomes
\(   \label{eq:E_PNPM}
    \mathcal{E}[D] = \int_{\Omega} \frac{\beta}{2} | \nabla D |^2 +  u \nabla \mu \cdot D \nabla \mu + \alpha |D|^\gamma \, \dx.
\) %beta diffusione dovuta al brownian motion
This functional describes the energy expenditure of a biological network
transporting charged particles. A generic example is represented by neural (brain) tissue
in animals and humans.
However, a quick inspection of equation \eqref{eq:gradient_flow_PNP}, or its
straightforward modification accounting for the presence of the metabolic term in \eqref{eq:E_PNPM},
reveals that, due to lack of a minimum principle, it does not guarantee preservation of positive (semi)definitness of the tensor $D$.
In a future work, we shall examine the well-posedness of the system in the case of small sources $S=S(x)$
and/or small time $t$.
Another option, inspired by \cite{AAFM, bookchapter, BHMR, HMP15, HMPS16}, is to make the ansatz
\(   \label{ansatz:m}
   D = r\mathbb{I} + m\otimes m,
\)
where $r=r(x) \geq r_0 >0$ is the background permeability of the medium
and the vector field $m=m(t,x)\in \mathbb{R}^d$ describes the local conductance of the network structure.
Note that $D$ taking the form \eqref{ansatz:m} has the eigenvalues $r(x) + |m|^2$ with eigenvector $m$,
and $r(x)$ with eigenvectors orthogonal to $m$.
Thus, it represents conduction along the direction $m$ with conductivity $r(x) + |m|^2$,
while the conduction in directions perpendicular to $m$ is due to the background permeability.

We, therefore, led to consider the following reformulation of \eqref{E_HuCai},
\(   %\label{eq:E_PNPM:m}
    \mathcal{E}[m] = \int_{\Omega} \frac{\beta}{2} |\grad m|^2 + r u  |\nabla \mu|^2 + u |m\cdot\nabla \mu|^2 + \frac{\alpha}{\gamma} |m|^{2\gamma} \, \dx, \notag
\)
with diffusivity $\beta>0$, metabolic constant $\alpha>0$ and metabolic exponent $\gamma>0$.
The functional \eqref{E_HuCai} is constrained by the drift-diffusion equation \eqref{eq:DD} with the diffusivity tensor $D$ given by \eqref{ansatz:m},
%\(  \label{eq:DD:m}
%          - \div\left[ r \nabla u + (m\cdot\nabla u) m + \z u r \nabla \varphi +  \z u  (m\cdot\nabla\varphi)m \right] = S,   % \qquad &\mathrm{in} \, \, \Omega, \\
%\)
and the Poisson equation \eqref{eq:Poisson} for the potential $\varphi$.

%%%%%%%%%%%%%%%%%%%%%%%%%%%%%%%%%%%%%%%%%%%%%%%
\section{General diffusive model} \label{sec:GDM}
We now derive the formal $L^2$-gradient flow of \eqref{funct0}
constrained by the elliptic problem \eqref{eq:ell}, \eqref{BC_u}.
Observe that multiplying \eqref{eq:ell} by $\Phi'(u)$ and integrating by parts one obtains
\[
   \int_\Omega \Phi''(u) \grad u \cdot D \grad u \,\d x = \int_\Omega S(x) \Phi'(u) \,\dx \ge 0,
\]
due to the convexity of $\Phi$ and since $\Phi'(c)=0$. Consequently, the functional $E[D]$ given by \eqref{funct0} can be written as
\begin{linenomath}
%\begin{subequations}
\begin{align}
\label{eq:energy_genmod}
    E[D] = \int_{\Omega} S(x) \Phi'(u) \, \dx.
\end{align}
%\end{subequations}
\end{linenomath}
%.

We then have the following result.

\begin{lemma}\label{lemma:gradient_flow_genmod}
The formal $L^2$-gradient flow of the energy functional \eqref{eq:energy_genmod} constrained by \eqref{eq:ell}, \eqref{BC_u} is given by
\begin{linenomath}
%\begin{subequations}
\begin{align}
    \label{eq:GF_diffusion}
    \part{D}{t} = \Phi''(u) \nabla u \otimes \nabla u + \frac{\nabla\sigma\otimes\nabla u + \nabla u\otimes\nabla\sigma}{2} \qquad \mbox{in } (0, \infty) \times \Omega, %\quad &\mathrm{in} \, \, \Omega,
\end{align}
%\end{subequations}
\end{linenomath}
%.
for the symmetric tensor valued diffusivity $D=D(t,x)$ and scalar $u=u(t,x)$,
with $\sigma=\sigma(t,x)$ the solution of the boundary value problem
\[   % \label{eq:auxiliary_problem_genmod}
    - \div (D \nabla \sigma) = \Phi'''(u) \nabla u \cdot D \nabla u \qquad \mathrm{in} \, \, (0, \infty) \times \Omega %& \quad &\mathrm{in} \, \, \Omega, \\
\]
subject to $\sigma = 0$ on $\partial \Omega$.
\end{lemma}

\begin{proof}
We expand $u=u^0 + \varepsilon u^1 + O(\varepsilon^2)$ and $D=D^0 + \varepsilon D^1 + O(\varepsilon^2)$ with $\varepsilon\in\R$,
where $D^0$ is a symmetric positive definite tensor and $D^1$ is symmetric.
Inserting into \eqref{eq:ell}, we obtain at zeroth-order
\(
      \label{eq:0order_system_genmod}
        - \diver(D^0 \nabla u^0 ) = S, %\quad &\mathrm{in} \, \, \Omega, \\
 \)
subject to $u^0 = c$ on $\partial \Omega$.
Collecting terms of first order in $\varepsilon$, we have
    \begin{linenomath}
    %\begin{subequations}
    \begin{align}
        \label{eq:1order_system_genmod}
        -\div(D^0 \nabla u^1 + D^1 \nabla u^0) = 0,
    \end{align}
    %\end{subequations}
    \end{linenomath}
subject to $u^1 = 0$ on $\partial \Omega$.
Multiplication of \eqref{eq:1order_system_genmod} by a sufficiently smooth function $v=v(x)$,
vanishing on the boundary $\partial\Omega$, and integration by parts leads to the useful identity
    \begin{linenomath}
    %\begin{subequations}
    \begin{align}
    \label{eq:identity1_genmod}
        \int_{\Omega} u^1 [-\diver (D^0 \nabla v)] \, \dx =\int_{\Omega} \nabla v \cdot D^0 \nabla u^1 \,\dx =  - \int_{\Omega}  \nabla v \cdot D^1 \nabla u^0 \,\dx.
    \end{align}
    %\end{subequations}
    \end{linenomath}

  Next, we calculate the first variation of $E$ in the direction $D^1$, % through a Taylor expansion, that is
    \begin{linenomath}
    %\begin{subequations}
    \begin{align}
        \frac{\delta E[D^0]}{\delta D}(D^1) &= \frac{d}{d\varepsilon} E[D^0 + \varepsilon D^1] \Bigr|_{\varepsilon = 0} \notag \\
        &= \frac{d}{d\varepsilon} \int_{\Omega} S(x) \Phi'(u^0 + \varepsilon u^1) \dx \Bigr|_{\varepsilon = 0}  \notag \\ &= \frac{d}{d\varepsilon} \int_{\Omega} S(x) \big[ \Phi'(u^0) + \varepsilon u^1 \Phi''(u^0) \big] \dx \Bigr|_{\varepsilon = 0} \notag \\ &= \int_{\Omega} u^1 S(x) \Phi''(u^0) \dx. \notag
    \end{align}
    %\end{subequations}
    \end{linenomath}
Using \eqref{eq:0order_system_genmod} and integrating by parts we obtain
    \begin{linenomath}
    %\begin{subequations}
    \begin{align}
    \label{eq:lemma1_step1_genmod}
         \frac{\delta E[D^0]}{\delta D}(D^1) &= \int_{\Omega} \nabla \big[ u^1 \Phi''(u^0) \big] \cdot D^0 \nabla u^0 \dx \notag \\ 
         &= \int_{\Omega} \nabla u^1 \cdot D^0 \nabla [\Phi'(u^0)] \dx + \int_{\Omega} u^1 \Phi'''(u^0) \nabla u^0 \cdot D^0  \nabla u^0 \dx =: I_1 + I_2. 
    \end{align}
    %\end{subequations}
    \end{linenomath}
    We apply \eqref{eq:identity1_genmod} with $v = \Phi'(u^0)$, noting that $\Phi'(u^0) =0$ on $\partial\Omega$ due to \eqref{BC_u} and \eqref{eq:equilibrium}, to evaluate
    \begin{linenomath}
    %\begin{subequations}
    \begin{align}
    \label{eq:I1_genmod}
        I_1 := \int_{\Omega} \nabla \Phi'(u^0) \cdot D^0 \nabla u^1 \dx = - \int_{\Omega} \nabla \Phi'(u^0) \cdot D^1 \nabla u^0 \,\dx. 
    \end{align}
    %\end{subequations}
    \end{linenomath}
    To evaluate $I_2$ we define $\sigma$ as the solution of the elliptic problem
    \begin{linenomath}
    %\begin{subequations}
    \begin{align}
       - \diver(D^0 \nabla \sigma) = \Phi'''(u^0) \nabla u^0 \cdot D^0 \nabla u^0, \notag
    \end{align}
    %\end{subequations}
    \end{linenomath}
    subject to homogeneous Dirichlet boundary condition on $\partial\Omega$.
Using \eqref{eq:identity1_genmod} with $v:=\sigma$ we arrive at
    \begin{linenomath}
    %\begin{subequations}
    \begin{align}
    \label{eq:I2_genmod}
        I_2 := \int_{\Omega} u^1 \Phi'''(u^0) \nabla u^0 \cdot D^0 \nabla u^0 \dx = - \int_{\Omega} \nabla \sigma \cdot D^1 \nabla u^0 \,\dx.
    \end{align}
    %\end{subequations}
    \end{linenomath}
    Due to the symmetry of $D^1$, we have
    \[
       I_2 &=& - \frac12 \int_{\Omega} \nabla \sigma \cdot D^1 \nabla u^0 + \nabla u^0 \cdot D^1 \nabla\sigma \,\dx \notag \\
          &=& - \int_{\Omega} D^1 : \frac{\nabla\sigma\otimes\nabla u^0 + \nabla u^0\otimes\nabla\sigma}{2} \,\dx,
    \]
        where the symbol $:$ denotes the contraction product of tensors, i.e., $A:B = \mathrm{tr}(A B^T)$. 
    Substituting \eqref{eq:I1_genmod} and \eqref{eq:I2_genmod} into \eqref{eq:lemma1_step1_genmod}, we finally get
    \begin{linenomath}
    %\begin{subequations}
    \begin{align}
        \frac{\delta E[D^0]}{\delta D}(D^1) &= - \int_{\Omega} D^1:\bigg[ \nabla u^0 \otimes \nabla \Phi'(u^0)
            + \frac{\nabla\sigma\otimes\nabla u^0 + \nabla u^0\otimes\nabla\sigma}{2} \bigg] \, \dx \notag \\
        &= - \int_{\Omega} D^1 : \bigg[ \Phi''(u^0) \nabla u^0 \otimes \nabla u^0 + \frac{\nabla\sigma\otimes\nabla u^0 + \nabla u^0\otimes\nabla\sigma}{2}\bigg] \, \dx. \notag
    \end{align} 
    %\end{subequations}
    \end{linenomath}
 \end{proof}

%%%%%%%%%%%%%%%%%%%%%%%%%%%%%%%%%%%%%%%%%%%%%%%%%%
\section{Drift-diffusion model}\label{sec:DD}
In this Section we extend the model to the drift-diffusion setting,
where the electrically charged particles are subject to a prescribed smooth stationary electric potential $\varphi=\varphi(x)$ with $\varphi = 0$ on $\partial \Omega$.
Then, the local mass conservation is of the form
\(  \label{eq:DD:DD}
        - \div( D \nabla u + \z u D \nabla \varphi ) = S \qquad \mbox{in } \Omega,  % \qquad &\mathrm{in} \, \, \Omega, \\
\)
where $\z\in\R$ denotes the valence (electric charge) of the particles.
Equation \eqref{eq:DD:DD} is subject to the Dirichlet boundary condition
\[  %   \label{eq:DD:DD:BC}
        u = c \qquad \mbox{on } \partial \Omega,
 \]
where the constant $c$ is an equilibrium, i.e., $\Phi'(c) = 0$.
%It is easily calculated that the dissipation of the Helmholtz entropy \eqref{eq:helmholtz_physmod}
%along the solutions of the parabolic version of the drift-diffusion equation \eqref{eq:DD:DD} with $S=0$ is of the form
%\[  \frac{d}{dt} H[u] = - \int_{\Omega} \Phi''(u) \nabla u \cdot \left(  D \nabla u + z u D \nabla \varphi \right) \,\dx.  \]
%Consequently,
%We introduce the entropy loss functional
%\(  \label{funct_dd}
%   E[D] := \int_{\Omega} \Phi''(u) \nabla u \cdot  D \left(\nabla u + \z u \nabla \varphi \right) \,\dx.
%\)
Defining the potential
\(  \label{w}
    w := e^{\z \varphi(x)} u,
\)
equation \eqref{eq:DD:DD} transforms into
\(  \label{eq:NP2_ddmod}
        - \div( e^{-\z\varphi} D \nabla w ) = S \qquad \mathrm{in} \, \, \Omega,
 \)
subject to $w = c$ on $\partial\Omega$.
%The entropy loss functional \eqref{funct_dd} then takes the form
We then consider the entropy loss functional
\(   \label{eq:energy_NP_ddmod}
    E[D] = \int_{\Omega} e^{-\z \varphi} \Phi''(w) \nabla w \cdot D \nabla w \dx. % = \int_{\Omega} S(x) \Phi'(w) \dx. \notag
\)

Following similar steps as in the proof of Lemma \ref{lemma:gradient_flow_genmod}, we derive
the $L^2$-gradient flow of \eqref{eq:NP2_ddmod}--\eqref{eq:energy_NP_ddmod}.

\begin{lemma}\label{gradient_flow_NP_ddmod}
The formal $L^2$-gradient flow of the functional \eqref{eq:energy_NP_ddmod}
constrained by \eqref{eq:NP2_ddmod} is given by 
\(  \label{eq:gradient_flow_NP_ddmod}
    \part{D}{t} = e^{-\z \varphi} \left[ \Phi''(w) \nabla w \otimes \nabla w + \frac{\nabla w \otimes \nabla \sigma + \nabla \sigma \otimes \nabla w}{2} \right] \qquad \mbox{in } (0, \infty) \times \Omega,
%    &= e^{z \varphi} \Phi''(e^{z \varphi} u) ( \nabla u + z u \nabla \varphi ) \otimes ( \nabla u + z u \nabla \varphi ) + (\nabla u + z u \nabla \varphi) \otimes \nabla \sigma,
\)
with $w$ given by \eqref{w} and $\sigma$ a solution of
\(   \label{eq:auxiliary_problem_NP_ddmod}
    - \div(D \nabla \sigma) + \z \nabla \varphi \cdot D \nabla \sigma = \Phi'''(w) \nabla w \cdot D \nabla w \qquad \mathrm{in} \, \, (0, \infty) \times \Omega % \quad &\mathrm{in} \, \, \Omega, \\
\)
subject to the homogeneous Dirichlet boundary condition $\sigma = 0$ on $\partial \Omega$.
\end{lemma}

\begin{proof}
See the proof of Lemma \ref{lemma:gradient_flow_genmod}.
\end{proof}

To examine convexity properties of the functional \eqref{eq:energy_NP_ddmod} constrained by \eqref{eq:NP2_ddmod},
we calculate its second-order variation.

\begin{lemma}\label{lem:conv}
The second-order variation of \eqref{eq:energy_NP_ddmod} coupled to \eqref{eq:NP2_ddmod} is given by
\(  \label{variation}
      \frac{\delta^2 E[D^0]}{\delta D^2}(D^1, D^1) = \int_\Omega e^{-\z\varphi} \grad\left( \Phi''(w^0) w^2 + \Phi'''(w^0) (w^1)^2 \right)\cdot D^0\grad w^0 \,\dx.
\)
where $w^0$ is a solution of
\[
   - \div( e^{-\z\varphi} D^0 \nabla w^0 ) = S \qquad \mbox{in } \Omega,
\]
$w^0 = c$ on $\partial\Omega$ and $w^1$, $w^2$ are defined by
\(    \label{c1}
   - \div \left[ e^{-\z\varphi} \left(D^1\grad w^0 + D^0 \grad w^1 \right)\right] = 0 \qquad \mathrm{in} \, \, \Omega, \\
     \label{c2}
   - \div \left[ e^{-\z\varphi} \left(D^1\grad w^1 + \frac12 D^0 \grad w^2 \right)\right] = 0 \qquad \mathrm{in} \, \, \Omega,
\)
subject to homogeneous Dirichlet boundary conditions.
\end{lemma}

\begin{proof}
We expand $D=D^0 + \varepsilon D^1 + O(\varepsilon^2)$ with $\varepsilon\in\R$,
where $D^0$ is a symmetric positive definite tensor and $D^1$ is symmetric.
By Taylor expansion we have
\[
   w[D^0 + \eps D^1] = w[D^0] + \eps \frac{\delta w[D^0]}{\delta D}(D^1) + \frac{\eps^2}{2} \frac{\delta^2 w[D^0]}{\delta D^2}(D^1, D^1) + O(\eps^3),
\]
and we denote
\[
   w^0 := w[D^0], \qquad w^1 := \frac{\delta w[D^0]}{\delta D}(D^1), \qquad w^2 := \frac{\delta^2 w[D^0]}{\delta D^2}(D^1, D^1).
\]
Collecting the $O(1)$ terms in \eqref{eq:NP2_ddmod} gives
\[
   - \div( e^{-\z\varphi} D^0 \nabla w^0 ) = S.
\]
The first-order terms give
\[
   - \div \left[ e^{-\z\varphi} \left(D^1\grad w^0 + D^0 \grad w^1 \right)\right] = 0,
\]
which is \eqref{c1}, and the second-order terms
\[
   - \div \left[ e^{-\z\varphi} \left(D^1\grad w^1 + \frac12 D^0 \grad w^2 \right)\right] = 0,
\]
which is \eqref{c2}.

Multiplication of \eqref{eq:NP2_ddmod} by $\Phi'(w)$ and integration by parts yields
\[
   E[D] = \int_\Omega S \Phi'(w) \d x,
\]
so that we have
\[
      \frac{\delta^2 E[D^0]}{\delta D^2}(D^1, D^1) = \totk{}{\eps}{2} \left. \int_\Omega S \Phi'\left( w[D^0 + \eps D^1] \right) \d x \right|_{\eps=0} . %\\
       %  &=& \totk{}{\eps}{2} \left. \int_\Omega S \Phi'\left( w^0 + \eps w^1 + \frac{\eps^2}{2} w^2 \right) \d x \right|_{\eps=0}.
\]
Taylor expansion gives
\[
   \Phi'\left( w[D^0 + \eps D^1] \right) &=& \Phi'\left( w^0 + \eps w^1 + \frac{\eps^2}{2} w^2 \right) + O(\eps^3) \\
     &=& \Phi'(w^0) + \eps \Phi''(w^0)w^1 + \frac{\eps^2}{2} \left( \Phi''(w^0) w^2 + \Phi'''(w^0) (w^1)^2 \right) + O(\eps^3).
\]
Consequently,
\[
      \frac{\delta^2 E[D^0]}{\delta D^2}(D^1, D^1) = \int_\Omega S \left( \Phi''(w^0) w^2 + \Phi'''(w^0) (w^1)^2 \right) \,\dx,
\]
and using \eqref{eq:NP2_ddmod} again, integration by parts results in
\[
      \frac{\delta^2 E[D^0]}{\delta D^2}(D^1, D^1) = \int_\Omega e^{-\z\varphi} \grad\left( \Phi''(w^0) w^2 + \Phi'''(w^0) (w^1)^2 \right)\cdot D^0\grad w^0 \,\dx.
\]
\end{proof}

We observe that for general (convex) $\Phi$ the result of Lemma \ref{lem:conv} does not directly imply convexity of $E[D]$.
However, for $\Phi(w) = w^2/2$ we have
\[
      \frac{\delta^2 E[D^0]}{\delta D^2}(D^1, D^1) = \int_\Omega e^{-\z\varphi} \grad w^2\cdot D^0\grad w^0 \,\dx.
\]
Multiplication of \eqref{c2} by $w^0$ and integration by parts gives
\[
     \frac12 \int_\Omega e^{-\z\varphi} \grad w^0 \cdot D^0 \grad w^2 \,\dx  =  - \int_\Omega e^{-\z\varphi} \grad w^0\cdot\D^1\grad w^1,
\]
and multiplication of \eqref{c1} by $w^1$ and integration by parts yields
\[
   \int_\Omega e^{-\z\varphi} \grad w^1\cdot\D^1\grad w^0 = - \int_\Omega e^{-\z\varphi} \grad w^1 \cdot D^0 \grad w^1 \,\dx.
\]
Consequently, recalling the symmetry and positive semidefinitness of $D^0$, we have
\[
      \frac{\delta^2 E[D^0]}{\delta D^2}(D^1, D^1) = \int_\Omega e^{-\z\varphi} \grad w^1 \cdot D^0 \grad w^1 \,\dx \geq 0.
\]
We conclude that for $\Phi(w) = w^2/2$ the minimization problem is convex.

For the following we denote by $S_d$ the space of symmetric real $d \times d$-matrices and, for $\alpha > 0$:
\[
    L^2_{\alpha} (\Omega; S_d) := \{ D \in L^2(\Omega; S_d) \, \, \mathrm{s. t.} \, \, D \ge \alpha I \, \, \mathrm{a.e.} \,  \mathrm{on} \, \, \Omega \}.
\]
We prove:
\begin{theorem}
Let $\varphi \in L^{\infty}(\Omega)$ and define the functional $E : L^2 (\Omega; S_d) \rightarrow (-\infty, +\infty] $ by
\[
    E[D] := \begin{cases}
    \int_{\Omega} e^{-\z \varphi} \nabla w \cdot D \nabla w \dx, \quad &\mathrm{if} \, \, D \in L^2_{\alpha} (\Omega; S_d), \\
    \infty \quad &\mathrm{otherwise},
    \end{cases}
\]
where $-\nabla \cdot (e^{-\z\varphi} D \nabla w ) = S$ in $\Omega$, $w=c$ on $\partial\Omega$, with $S \in L^2(\Omega)$. Then, for each $D_I \in L^2_{\alpha} (\Omega; S_d)$, the gradient flow of $E$, given by
\(
     \label{eq:GF_D1}
    \part{D}{t} = e^{-\z \varphi} \nabla w \otimes \nabla w \qquad &\mbox{on }& (0, \infty) \times \Omega, \\
    D(0, x) = D_I (x) \qquad &\mbox{on }& \Omega,
     \label{eq:GF_D2}
\)
coupled to \eqref{eq:NP2_ddmod}, exists for all $t>0$.

\begin{proof}
For $S \in L^2_{\alpha}(\Omega; S_d)$ an application of the Lax-Milgram Lemma in the space 
\[
    H := \left\{ v \in H^1_0 (\Omega) \, \, \mathrm{s.t.} \, \, \int_{\Omega} e^{-\z \varphi} \nabla v \cdot D \nabla v \, \dx < \infty \right\}
\] 
shows the existence of a unique solution $w \in H + c$ of the equation $-\nabla \cdot (e^{-\z \varphi} D \nabla w) = S$ in $\Omega$, $w = c$ on $\partial\Omega$. Thus, $E$ is well defined. Moreover, it is easy to prove that $E$ is $C^2$ on its domain of definition.

Obviously, $L^2_{\alpha}(\Omega; S_d)$ is a closed and convex subset of $L^2(\Omega; S_d)$ and E is convex and lower semicontinuous on $L^2(\Omega; S_d)$ (see the proof of Proposition 3.3 in \cite{CMS-2022} for the lower semicontinuity proof). The statement then follows from standard theory, see e.g. \cite{ABS}.
\end{proof}
\end{theorem}

Let us note that the right-hand side in \eqref{eq:GF_D1} is a positive semidefinite matrix.
Consequently, during the evolution induced by \eqref{eq:GF_D1} the solution $D=D(t)$
'stays away' from the boundary of the set $L^2_{\alpha} (\Omega; S_d)$, which consists of
positive semidefinite $L^2$-integrable tensors. In particular, if $D_I \geq \alpha I$ almost everywhere
on a subset $U$ of $\Omega$ of positive Lebesque measure, then $D(t,x) \geq \alpha I$
almost everywhere on $U$ for all $t\geq 0$.

Note that the gradient flow equation \eqref{eq:GF_D1} for $D=D(t)$, coupled to \eqref{eq:NP2_ddmod},
does not become stationary (except for the trivial case $S\equiv 0$, when $\grad w \equiv 0$).
%as $D \rightarrow \infty$, instead $D$ is pointwise increasing in time (in the sense of positive-definite matrices).
This is mitigated by adding a relaxation term for $D$ to the right-hand side of \eqref{eq:GF_D1}.
Following the suggestion of \cite{HMP15, HMPS16, CMS-2022, Hu-Cai}, we choose the power-law $|D|^{\gamma-2} D$ with $\gamma> 1$.
% and possibly a diffusion term as in \eqref{E_HuCai}.
%In particular, the choice $\gamma > 1$, %$\beta=0$ (vanishing diffusion),
%and, for simplicity, $\varphi \equiv 0$ % and $\alpha = 1$ (for the sake of simplicity) we obtain the stationary state equation from \eqref{eq:GF_Dir_complete}
Setting  $\varphi \equiv 0$ for simplicity, we arrive at the following stationary version of \eqref{eq:GF_D1},
\[
    \grad w \otimes \grad w = |D|^{\gamma-2} D.
\]
Inserting into \eqref{eq:NP2_ddmod}, we arrive at
\[
    - \diver ( |\grad{w}|^{\frac{2}{\gamma-1}} \nabla w ) = S \qquad \mbox{in } \Omega,
\]
%(see \cite{HMP15, HMPS16, CMS-2022}).
which is the p-Laplace equation with $p=\frac{2\gamma}{\gamma-1} > 1$.
In this context let us point out the important works \cite{Vazquez1, Vazquez2} of Juan Luis V{\'a}zquez.

%As anticipated in the introduction and in the second-order variation calculation in Lemma 3, we can neither prove that the energy functional is convex nor $D$ remains non-negative for general entropy dissipation.
A significant problem for proving well-posedness of the system \eqref{eq:NP2_ddmod}, \eqref{eq:gradient_flow_NP_ddmod}, \eqref{eq:auxiliary_problem_NP_ddmod}
with general convex entropy loss densities $\Phi=\Phi(w)$
is the fact that we are not able to establish preservation of nonnegativity of the tensor $D=D(t)$.
However, modeling considerations \cite{AAFM, bookchapter, BHMR, HMP15, HMPS16} motivate us to make
the ansatz \eqref{ansatz:m} for $D$, namely
\[
   D = r\mathbb{I} + m\otimes m,
\]
with the regularization parameter $r=r(x) \geq r_0 >0$ (background permeability of the medium)
and the vector field $m=m(t,x)\in \mathbb{R}^d$ (local conductance of the network structure).
Then, the Fokker-Planck equation \eqref{eq:NP2_ddmod} transforms into
\( \label{eq:NP2m_ddmod}
     - \div \big[ e^{-\z\varphi} (r I + m \otimes m) \nabla w \big] = S \qquad \mathrm{in} \, \, \Omega,
\)
subject to $w = c$ in $\partial\Omega$. Similarly, we recast the entropy loss functional \eqref{eq:energy_NP_ddmod} as
\( \label{eq:energy_NPm_ddmod}
    E[m] = \int_{\Omega} e^{-\z \varphi} \Phi''(w) \grad w \cdot (rI + m \otimes m) \grad w \,\dx.
\)
Note that a multiplication of $\eqref{eq:NP2m_ddmod}$ by $\Phi'(w)$ and an integration by parts gives
\[
   E[m] = \int_{\Omega} S(x) \Phi'(w) \dx.
\]
We have the following form for the $L^2$-gradient flow of the system \eqref{eq:NP2m_ddmod} -- \eqref{eq:energy_NPm_ddmod}.

\begin{lemma}
The formal $L^2$-gradient flow of the functional \eqref{eq:energy_NPm_ddmod} constrained by \eqref{eq:NP2m_ddmod} is given by
\[ \label{eq:GF_FP_m}
    \frac{\partial m}{\partial t} = e^{-\z \varphi} \big[ 2 \Phi''(w) (\grad w \cdot m) \grad w + (\grad \sigma \cdot m) \grad w + (\grad w \cdot m) \grad \sigma \big] \qquad \mathrm{in} \, \, (0, \infty) \times \Omega,
\]
with $w$ given by \eqref{w} and $\sigma$ a solution of 
\[
    - \grad \cdot \big[ (r I + m \otimes m) \grad \sigma \big] + \z \grad \varphi \cdot (rI + m \otimes m) \grad \sigma = \Phi'''(w) \grad w \cdot (rI + m \otimes m) \grad w % \, \, \mathrm{in} \, (0, \infty) \times \Omega
\]
subject to the homogeneous Dirichlet boundary condition $\sigma = 0$ on $\partial \Omega$.
\begin{proof}
We only sketch the proof here, following the lines of the proof of Lemma \ref{lemma:gradient_flow_genmod},
where we substitute for $u^0= e^{-\z \varphi} w^0$, $u^1= e^{-\z \varphi} w^1$, $D^0 = rI + m^0 \otimes m^0$ and $D^1 = m^0 \otimes m^1 + m^1 \otimes m^0$.
% obtained after the usual expansion in $\varepsilon$ of $m$.

Combining \eqref{eq:lemma1_step1_genmod}, \eqref{eq:I1_genmod} and \eqref{eq:I2_genmod} and setting $\sigma$ as a solution of 
\[ 
    - \grad \cdot \big[ e^{-\z \varphi} (rI + m^0 \otimes m^0) \grad \sigma \big] = e^{-\z \varphi} \Phi'''(w^0) \grad w^0 \cdot (rI + m^0 \otimes m^0) \grad w^0,
\]
we get
\[
    \frac{\delta E[m^0]}{\delta m}(m^1) &=& - \int_{\Omega} e^{- \z \varphi} \grad \Phi'(w^0) \cdot (m^0 \otimes m^1 + m^1 \otimes m^0) \grad w^0 \dx \\
    &-& \int_{\Omega} e^{-\z \varphi} \grad \sigma \cdot (m^0 \otimes m^1 + m^1 \otimes m^0) \grad w^0 \dx \\
    &=& - \int_{\Omega} m^1 \cdot e^{-\z \varphi} \big[ 2 \Phi''(w^0) ( \grad w^0 \cdot m^0) \grad w^0 + (\grad \sigma \cdot m) \grad w + (\grad w \cdot m) \grad \sigma \big] \dx.
\]
\end{proof}
\end{lemma}

Let us now examine the (non)convexity of the functional \eqref{eq:energy_NPm_ddmod}
with $\Phi(w)=w^2/2$.

\begin{lemma}
Denote $\d\zeta := e^{-\z\varphi} \dx$.
The second-order variation of the energy
\(   \label{Em}
   E[m] := \int_\Omega  r|\grad w|^2 + |m\cdot\grad w|^2 \,\d\zeta
\)
in direction $m^1\in H^1_0(\Omega)$, where $w=w[m]$ is the solution of \eqref{eq:NP2m_ddmod}, reads
\[
   \frac{\delta^2 E[m^0]}{\delta m^2}(m^1, m^1) = 2 \int_\Omega r \left| \grad \frac{\delta w[m^0]}{\delta m}(m^1) \right|^2
       + \left| m^0\cdot\grad \frac{\delta w[m^0]}{\delta m}(m^1) \right|^2 
       - \left| m^1\cdot\grad w \right|^2 \,\d\zeta.
\]
\end{lemma}

\begin{proof}
Using $w$ as test function in the weak formulation of \eqref{eq:NP2m_ddmod} gives
\[
   E[m] = \int_\Omega S w \d x,
\]
so that
\[
   \frac{\delta^2 E[m^0]}{\delta m^2}(m^1, m^1) = \int_\Omega S \frac{\delta^2 w[m^0]}{\delta m^2}(m^1, m^1).
\]
Let us denote
\[
    w^0 := w[m], \qquad w^1 := \frac{\delta w[m^0]}{\delta m}(m^1), \qquad w^2 := \frac{\delta^2 w[m^0]}{\delta m^2}(m^1,m^1).
\]
We use $w^1$ as a test function in the weak formulation of \eqref{eq:NP2m_ddmod}
and calculate the first-order variation in direction $m^1 \in H^1_0(\Omega)$, which leads to
\(
    \nonumber
   \int_\Omega S w^2 \,\dx &=&
       \int_\Omega r \left| \grad w^1 \right|^2 + r \grad w^0 \cdot\grad w^2 \,\d\zeta \\
       \label{eq:Sw2}
       &+& \int_\Omega (m^1\cdot\grad w^0)\left( m^0 \cdot\grad w^1 \right)
             + \left| m^0\cdot\grad w^1 \right|^2 \\
            && \quad
              + (m^0\cdot\grad w^0)\left( m^1 \cdot \grad w^1 \right)
                + (m^0\cdot\grad w^0)\left( m^0\cdot \grad w^2 \right) \,\d\zeta.
                \nonumber
\)
The first-order variation of the weak formulation of \eqref{eq:NP2m_ddmod} with test function $\xi\in H^1_0(\Omega)$ reads
\(  \label{eq:varPoisson}
   \int_\Omega r \grad w^1\cdot\grad\xi &+& (m^1\cdot\grad w^0)(m^0\cdot\grad\xi) + (m^0\cdot\grad w^1)(m^0\cdot\grad\xi) \notag \\ &+& (m^0\cdot\grad w^0)(m^1\cdot\grad\xi) \,\d \zeta = 0,
\)
and setting $\xi:=w^1$ gives
\( \label{eq:varPoisson_w1}
   \int_\Omega r |\grad w^1|^2 + (m^1\cdot\grad w^0)(m^0\cdot\grad w^1) + \left|m^0\cdot\grad w^1 \right|^2 + (m^0\cdot\grad w^0)(m^1\cdot\grad w^1) \,\d \zeta = 0.
\)
Inserting into \eqref{eq:Sw2} gives
\[
   \int_\Omega S w^2 \,\dx =
       \int_\Omega r \grad w^0 \cdot\grad w^2 + (m^0\cdot\grad w^0)\left( m^0\cdot \grad w^2 \right) \,\d\zeta.
\]
Now we again take a variation of \eqref{eq:varPoisson} in direction $m^1$ and use $\xi:=w^0$ as the test function,
\[
   \int_\Omega r \grad w^2\cdot\grad w^0 +  (m^0\cdot\grad w^0)\left( m^0\cdot \grad w^2 \right) \,\d\zeta = \\
    = - 2 \int_\Omega (m^1\cdot\grad w^1)(m^0\cdot\grad w^0) + \left| m^1\cdot\grad w^0 \right|^2 + (m^0\cdot\grad w^1)(m^1\cdot\grad w^0) \,\d\zeta.
\]
Consequently,
\[
   \int_\Omega S w^2 \,\dx = - 2 \int_\Omega (m^1\cdot\grad w^1)(m^0\cdot\grad w^0) + \left| m^1\cdot\grad w^0 \right|^2 + (m^0\cdot\grad w^1)(m^1\cdot\grad w^0) \,\d\zeta.
\]
Using \eqref{eq:varPoisson_w1}, we finally arrive at
\[
    \int_\Omega S w^2 \,\dx = 2 \int_\Omega r \left| \grad w^1 \right|^2 + \left| m^0\cdot\grad w^1 \right|^2 - \left| m^1\cdot\grad w^0 \right|^2 \,\d\zeta.
\]
\end{proof}

To gain a better insight into the convexity properties of the functional \eqref{Em},
let us recall the spatially one-dimensional case considered in \cite[Remark 2.4]{AAFM}.
We set $\Omega:=(0,1)$ and, for simplicity, $r(x) :\equiv 1$ and $\varphi := 0$.
Moreover, to expedite the calculation, we impose the mixed boundary conditions for $w$,
$\part{w}{x}(0) = w(1) = 0$. Then an integration of the Poisson equation \eqref{eq:NP2m_ddmod}
gives
\[
   \part{w^0}{x}(x) = - \frac{B(x)}{1+m^2},
\]
with $B(x) := \int_0^x S(\xi) \d\xi$. The Frechet derivative in direction $m^1$ reads then
\[
    \part{w^1}{x} = \part{}{x} \frac{\delta w[m^0]}{\delta m}(m^1) = \frac{2Bm^0}{(1+(m^0)^2)^2} m^1.
\]
Consequently, we have
\[
   \frac{\delta^2 E[m^0]}{\delta m^2}(m^1, m^1) 
      &=& 2 \int_\Omega  \left| \part{w^1}{x} \right|^2 + \left| m^0 \part{w^1}{x} \right|^2 - \left| m^1 \part{w^0}{x} \right|^2  \\
      &=& 2 \int_\Omega \frac{(m^1)^2 B^2}{(1+(m^0)^2)^3} \left( 3(m^0)^2 - 1 \right) \, \dx.
\]
Clearly, the sign of the second-order variation of $E[m^0]$ in any direction $m^1$ depends on $m^0$,
i.e., if $|m^0| \geq 1/\sqrt{3}$ then is non-negative; otherwise it is negative.
Therefore, $E$ is not convex on $L^2(\Omega)$.

%We observe that the sign of the second-order variation  does not depend on the choice of the direction $m^1$.
%Clearly, the  the convexity of $E[m^0]$ depends on the size of $(m^0)^2$; if $(m^0)^2 \geq 1/3$ almost everywhere in $(0,1)$, then $E[m^0]$ is convex.

%%%%%%%%%%%%%%%%%%%%%%%%%%%%%%%%%%%%%%%%%%%%%%%
\section{Poisson-Nernst-Planck model}\label{sec:PNP}
In this section we consider the convection-diffusion equation \eqref{eq:DD:DD} for an ion charge density $u$ with drift induced by the electrostatic field of the charged particles.
Consequently, the electric potential $\varphi=\varphi(t,x)$ is a solution of the Poisson equation
\(  \label{eq:Poisson_PNP}
   -\laplace\varphi = \z u \qquad \mbox{in } \Omega,
\)
where $\z\in\R$ is the particle charge. We prescribe homogeneous Dirichlet boundary condition for $\varphi$,
\(    \label{eq:dirichlet_varphi_ddmod}
    \varphi = 0 \qquad \mbox{on } (0,T] \times \partial \Omega.
\)
%The coupled system \eqref{eq:DD}, \eqref{eq:Poisson_PNP} then constitutes a Poisson-Nernst-Planck system.
As argued in, e.g., \cite{lu2011}, the entropy generator $\Phi=\Phi(u)$ for the Poisson-Nernst-Planck system
\eqref{eq:DD:DD}, \eqref{eq:Poisson_PNP} is given by $\Phi(u) = u (\ln u -1)$,
and the Helmholtz free energy takes the form \eqref{eq:helmholtz_ddmod}.
Moreover, observe that we have the equilibrium state $c=1$, i.e., $\Phi'(1) = 0$.

Introducing the quasi-Fermi energy level $\mu$ defined by \eqref{eq:mu_ddmod}, the parabolic version of the Poisson-Nernst-Planck system reads:
\(
   \label{eq:PNP_F1}
    \part{u}{t} - \div ( u D \nabla \mu) &=& 0 \qquad \mbox{in } (0, \infty) \times \Omega,\\
   \label{eq:PNP_F2}    
    - \Delta \varphi &=& \z u \qquad \mbox{in } (0, \infty) \times \Omega.
\)
It is a known result, see, e.g., \cite{BDP}, that the loss of the Helmholtz free energy \eqref{eq:helmholtz_ddmod} along the solutions
of the Poisson-Nernst-Planck system \eqref{eq:PNP_F1}--\eqref{eq:PNP_F2} is given by the functional
\(
    \label{eq:PNP_energy_functional1_ddmod}
    \mathcal{E}[D] = \int_{\Omega} u \nabla \mu \cdot D \nabla \mu \, \dx.
\)
For the convenience of the reader, we detail the calculation here.

\begin{lemma}\label{lemma:dissipation_ddmod}
We have, along the solutions of \eqref{eq:PNP_F1}--\eqref{eq:PNP_F2},
\[
    \frac{d}{dt} \mathcal{H}(u, \varphi) = - \mathcal{E}[D],
\]
where $\mathcal{H}(u, \varphi)$ is the Helmholtz free energy \eqref{eq:helmholtz_ddmod}
and $\mathcal{E}[D]$ is given by \eqref{eq:PNP_energy_functional1_ddmod}.
\end{lemma}

\begin{proof}
We calculate
\[
     \frac{d}{dt} \mathcal{H}(u, \varphi) &=& \int_{\Omega} \ln(u) \part{u}{t} + \nabla \varphi \cdot \nabla \part{\varphi}{t} \,\dx \\
       &=& \int_{\Omega} \ln(u)\part{u}{t}  - \varphi\laplace  \part{\varphi}{t}\, \dx \\
       &=& \int_{\Omega} \left( \ln(u) + \z \varphi \right) \part{u}{t} \, \dx,
\]
where we used \eqref{eq:PNP_F2} and integrated by parts.
Substituting for $\part{u}{t}$ from \eqref{eq:PNP_F1} and integrating by parts again gives
\[
     \frac{d}{dt} \mathcal{H}(u, \varphi) &=& - \int_{\Omega} u \grad \left( \ln(u) + \z \varphi \right)\cdot D\grad\mu \, \dx \\
                 &=& - \int_{\Omega}  \left( \grad u + \z u \grad\varphi \right) \cdot D\grad\mu \, \dx \\
                 &=&  - \int_{\Omega} u \nabla \mu \cdot D \nabla \mu \, \dx.
\]
Note that the boundary terms in the partial integration steps vanish due to the
homogeneous boundary condition \eqref{eq:dirichlet_varphi_ddmod} for $\varphi$.
\end{proof}

We now derive the $L^2$-gradient flow of the loss functional \eqref{eq:PNP_energy_functional1_ddmod}
coupled to the stationary Poisson-Nernst-Planck system formed by \eqref{eq:Poisson_PNP} and
\(
    \label{eq:PNP_elliptic1}
    - \div(u D  \nabla \mu) = S \qquad &\mathrm{in} \, \, \Omega,
\)
subject to the boundary conditions
\(
   \label{eq:PNP_bc_ddmod}
    u = 1, \quad \varphi= 0, \qquad \mathrm{on} \, \, \partial \Omega.
\)
Note that \eqref{eq:PNP_bc_ddmod} implies $\mu = 0$ on $\partial \Omega$.
Consequently, multiplication of \eqref{eq:PNP_elliptic1} by $\mu$ and integration by parts gives
\(
   \label{eq:PNP_energy_functional2_ddmod}
    \mathcal{E}[D] = \int_{\Omega} u \nabla \mu \cdot D \nabla \mu \dx = \int_{\Omega} S(x) \mu  \, \dx.
\)

\begin{lemma}\label{GF_PNP}
The formal $L^2$-gradient flow of the functional \eqref{eq:PNP_energy_functional2_ddmod} constrained by \eqref{eq:PNP_elliptic1}--\eqref{eq:PNP_bc_ddmod} is given by
\begin{linenomath}
%\begin{subequations}
\begin{align}
\label{eq:gradient_flow_PNP_ddmod}
    \part{D}{t} = u \nabla \mu \otimes \nabla \mu + u \frac{\nabla\mu \otimes \nabla\sigma + \nabla\sigma \otimes \nabla\mu}{2} \qquad &\mathrm{in} \, \, (0, \infty) \times \Omega,
\end{align}
%\end{subequations}
\end{linenomath}
coupled to the system for the auxiliary quantities $\sigma=\sigma(t,x)$ and $\eta=\eta(t,x)$,
\begin{linenomath}
%\begin{subequations}
\begin{align}
   %\label{eq:auxiliary_problem_PNP1}
    - \diver (D \nabla \sigma) + \z \nabla \varphi \cdot D \nabla \sigma - \z^2 \eta = \nabla \mu \cdot D \nabla \mu \qquad &\mathrm{in} \, \, (0, \infty) \times \Omega, \notag \\
   %\label{eq:auxiliary_problem_PNP2}
    - \Delta \eta = \diver(u D \nabla \sigma) \qquad &\mathrm{in} \, \, (0, \infty) \times \Omega \notag
\end{align}
%\end{subequations}
\end{linenomath}
subject to the boundary conditions
\(
   %\label{eq:auxiliary_problem_PNP_BC}
    \sigma = 0, \quad \eta=0 \qquad \mbox{on }  \partial \Omega. \notag
\)
\end{lemma}

\begin{proof}
Let us expand $D = D^0 + \varepsilon D^1 + O(\varepsilon^2)$,
where $D^0$ is a symmetric positive definite tensor and $D^1$ is symmetric.
Similarly, we expand the other relevant quantities $u$, $\varphi$ and $\mu$
 in terms of $\eps>0$.
% $u = u^0 + \varepsilon u^1 + O(\varepsilon^2)$,
%\varphi = \varphi^0 + \varepsilon \varphi^1 + o(\varepsilon^2)$
%$\mu = \mu^0 + \varepsilon \mu^1 + o(\varepsilon^2)$.
With \eqref{eq:mu_ddmod} we have $\mu^1 = \frac{u^1}{u^0} + \z \, \varphi^1$.
Collecting the zero-order terms in \eqref{eq:PNP_elliptic1} and \eqref{eq:Poisson_PNP}, we have
\[
     - \div( u^0 D^0 \nabla \mu^0) &=& S, \\
     - \Delta \varphi^0 &=& \z u^0.
\]
subject to the boundary conditions $u^0=1$ and $\varphi^0=0$ on $\partial\Omega$.
At first order in $\eps$ we obtain the system
\(
    \label{eq:1order1_system_PNP1}
        - \div( u^0 D^1 \nabla \mu^0 + u^1 D^0 \nabla \mu^0 + u^0 D^0 \nabla \mu^1) &=& 0  \\
    \label{eq:1order1_system_PNP2}
        - \Delta \varphi^1 &=& \z u^1,
 \)
subject to $u^1 = 0$ and $\varphi^1 = 0$ on $\partial \Omega$.
Note that \eqref{eq:1order1_system_PNP1} can also be rewritten in the form
\(       \label{eq:1order1_system_PNP1alt}
    - \div\big[ D^0 ( \nabla u^1 + \z u^1 \nabla \varphi^0 + \z u^0 \nabla \varphi^1) + D^1 ( \nabla u^0 + \z u^0 \varphi^0 )  \big] = 0 \qquad \mathrm{in} \, \, \Omega.
\)
%This system is uniquely solvable in $u^1$ for a given second-order, symmetric, positive semi-definite tensor $D^1$. Therefore, we introduce the corresponding bijective solution operator $\mathcal{A} : \R^{2 \times 2} \rightarrow \R $ such that $u^1 = \mathcal{A}(D^1)$.

Next, we calculate the first variation of $\mathcal{E}$ given by \eqref{eq:PNP_energy_functional2_ddmod} in the direction $D^1$,
\(
    \label{eq:1_lemmaPNP_ddmod}
    \frac{\delta \mathcal{E}[D^0]}{\delta D} (D^1) &=& \frac{d}{d\varepsilon} \int_{\Omega} S(x) (\mu^0 + \varepsilon \mu^1) \, dx \Bigr|_{\varepsilon = 0} \notag \\
    &=& \int_{\Omega} - \diver (u^0 D^0 \nabla \mu^0) \mu^1 \, \dx \notag \\
    &=& \int_{\Omega} u^0 \nabla \mu^1 \cdot D^0 \nabla \mu^0 \, \dx.
\)
Multiplication of the first-order system \eqref{eq:1order1_system_PNP1} by $\mu^0$ and integration by parts, recalling the symmetry of $D$, leads to
\[
    \int_{\Omega} u^0 \nabla \mu^1 \cdot D^0 \nabla \mu^0 \dx = - \int_{\Omega} u^0 \nabla \mu^0 \cdot D^1 \nabla \mu^0 + u^1 \nabla \mu^0 \cdot D^0 \nabla \mu^0 \dx.
\]
Substitution of the above identity into \eqref{eq:1_lemmaPNP_ddmod} gives
\(  \label{eq:2_lemmaPNP_ddmod}
    \frac{\delta \mathcal{E}[D^0]}{\delta D} (D^1) = - \int_{\Omega} u^0 \nabla \mu^0 \cdot D^1 \nabla \mu^0 + u^1 \nabla \mu^0 \cdot D^0 \nabla \mu^0 \dx.
\)
To evaluate the second term, we need to find a mapping between $u^1$ and $D^1$. For this sake, we multiply \eqref{eq:1order1_system_PNP1alt}
by a function $\sigma$ vanishing at the boundary $\partial\Omega$ and integrate by parts,
\[
     \int_{\Omega} \nabla \sigma \cdot \big[D^0 (\nabla u^1 + \z u^1 \nabla \varphi^0 + \z u^0 \nabla \varphi^1 ) \big] \, \dx &=&
         - \int_{\Omega} \nabla \sigma \cdot D^1 (\nabla u^0 + \z u^0 \nabla \varphi^0) \, \dx \\
    &=& - \int_{\Omega} u^0 \nabla \sigma \cdot D^1 \nabla \mu^0 \, \dx.
\]
After further integration by parts on the left-hand side, we obtain
\(    \label{eq:3_lemmaPNP_ddmod}
    - \int_{\Omega} u^1 \big[ \diver(D^0 \nabla \sigma) - \z \nabla \sigma \cdot D^0 \nabla \varphi^0 \big] + \z \varphi^1 \diver(u^0 D^0 \nabla \sigma) \dx
        = - \int_{\Omega} u^0 \nabla \sigma \cdot D^1 \nabla \mu^0 \dx.
\)
With the Poisson equation \eqref{eq:1order1_system_PNP2} we rewrite the third term of the left-hand side as
\[
    \int_{\Omega} \z \varphi^1 \diver(u^0 D^0 \nabla \sigma) \,\dx &=& - \int_{\Omega} \z^2 \Delta^{-1} \big[ u^1 \diver (u^0 D^0 \nabla \sigma) \big] \,\dx \\
    &=& - \int_{\Omega} \z^2 u^1 \Delta^{-1} \big[ \diver (u^0 D^0 \nabla \sigma) \big] \, \dx.
\]
Consequently, defining the function $\psi=\psi[\sigma]$,
\[
     \psi[\sigma] := - \diver(D^0 \nabla \sigma) + \z \nabla \sigma \cdot D^0 \nabla \varphi^0 + \z^2 \Delta^{-1} \diver( u^0 D^0 \nabla \sigma),
\] 
equation \eqref{eq:3_lemmaPNP_ddmod} becomes
%\(    \label{eq:4_lemmaPNP_ddmod}
%    \int_{\Omega} u^1 \big[-\diver(D^0 \nabla \sigma) + \z \nabla \sigma \cdot D^0 \nabla \varphi^0 + \z^2 \Delta^{-1} \diver(u^0 D^0 \nabla \sigma) \big] \,\dx
%        = - \int_{\Omega} u^0 \nabla \sigma \cdot D^1 \nabla \mu^0 \,\dx.
%\)
\[
    \int_{\Omega} u^1 \psi[\sigma] \, \dx = - \int_{\Omega} u^0 \nabla \sigma \cdot D^1 \nabla \mu^0 \, \dx.
\]
Hence, setting $\psi[\sigma] = \nabla \mu^0 \cdot D^0 \nabla \mu^0$,
% - which is the auxiliary System \eqref{eq:auxiliary_problem_PNP_ddmod} - 
we obtain
\[
     - \int_{\Omega} u^1 \nabla \mu^0 \cdot D^0 \nabla \mu^0 \,\dx = - \int_{\Omega} u^0 \nabla \sigma \cdot D^1 \nabla \mu^0 \,\dx.
\]
Substitution of the above identity into \eqref{eq:2_lemmaPNP_ddmod} leads to
\[
    \frac{\delta \mathcal{E}[D^0]}{\delta D} (D^1) = - \int_{\Omega} u^0 \nabla \mu^0 \cdot D^1 \nabla \mu^0 + u^0 \nabla \sigma \cdot D^1 \nabla \mu^0 \,\dx.
\]
Due to the symmetry of $D^1$, we have
\[
    \nabla \sigma \cdot D^1 \nabla \mu^0 &=& \frac{\nabla \sigma \cdot D^1 \nabla \mu^0 + \nabla \mu^0 \cdot D^1 \nabla \sigma}{2} \\
      &=& D^1 : \frac{\nabla \sigma \otimes \nabla \mu^0 + \nabla \mu^0 \otimes \nabla \sigma}{2}.
\]
We thus finally arrive at
\[
       \frac{\delta \mathcal{E}[D^0]}{\delta D} (D^1) =
          - \int_{\Omega} D^1 : \left[ u^0 \nabla \mu^0 \otimes \nabla \mu^0 + u^0 \frac{\nabla \sigma \otimes \nabla \mu^0 + \nabla \mu^0 \otimes \nabla \sigma}{2} \right] \,\dx,
\]
which directly gives \eqref{eq:gradient_flow_PNP_ddmod}.
\end{proof}
Once again, we cannot guarantee the non-negativity of $D$ for every time $t>0$ from Equation \eqref{eq:gradient_flow_PNP_ddmod}, so we employ the ansatz \eqref{ansatz:m}. Therefore, we obtain the following:
\begin{lemma}
The formal $L^2$-gradient flow of the functional \eqref{eq:PNP_energy_functional2_ddmod} constrained by \eqref{eq:PNP_elliptic1}--\eqref{eq:PNP_bc_ddmod} is given by
\begin{linenomath}
%\begin{subequations}
\begin{align}
\label{eq:gradient_flow_PNP_ddmod}
    \part{m}{t} = 2 u (m \cdot \grad \mu) \grad \mu + (m \cdot \nabla \mu) \nabla \sigma + (m \cdot \nabla \sigma) \nabla \mu \qquad &\mathrm{in} \, \, (0, \infty) \times \Omega,
\end{align}
%\end{subequations}
\end{linenomath}
coupled to the system for the auxiliary quantities $\sigma=\sigma(t,x)$ and $\eta=\eta(t,x)$,
\begin{linenomath}
%\begin{subequations}
\begin{align}
%\label{eq:auxiliary_problem_PNP1}
    - \diver \big[ (rI + m \otimes m) \nabla \sigma \big] + \z \nabla \varphi \cdot (rI + m \otimes m) \nabla \sigma - \z^2 \eta = \nabla \mu \cdot (rI + m \otimes m) \nabla \mu, \notag \\ %\qquad &\mathrm{in} \, \, (0, \infty) \times \Omega, \notag \\
   %\label{eq:auxiliary_problem_PNP2}
    - \Delta \eta = \diver \big[ u (rI + m \otimes m) \nabla \sigma \big], \notag %\qquad &\mathrm{in} \, \, (0, \infty) \times \Omega \notag
\end{align}
%\end{subequations}
\end{linenomath}
in $(0, \infty) \times \Omega$ and subject to the boundary conditions
\(
   %\label{eq:auxiliary_problem_PNP_BC}
    \sigma = 0, \quad \eta=0 \qquad \mbox{on }  \partial \Omega. \notag
\)
\begin{proof}
See the proof of Lemma \ref{GF_PNP}.
\end{proof}
\end{lemma}

%%%%%%%%%%%%%%%%%%%%%%%%%%%%%%%%%%%%%%%%%%%%%%%%%%%%%%
\section{Conclusions}

In this paper we systematically developed a class of self-regulating processes governed by the minimization of an entropy dissipation coupled to the conservation law for a quantity representing concentration of a chemical species, ions, nutrients or material pressure. The corresponding constrained $L^2$-gradient flows provide systems of partial differential equations that describe formation and evolution of biological transportation networks. We started with the derivation of the $L^2$-gradient flow for general entropy dissipations associated with a purely diffusive model. Further, we introduced the Fokker-Planck equation to describe charged ions, prescribing a given stationary electric potential, and computed the $L^2$ gradient flow for the drift-diffusion equation. We showed global existence and uniqueness of a solution in the case of a quadratic entropy generator due to the convexity and lower semicontinuity of the energy functional. Finally, we extended the model to the Poisson-Nernst-Planck (PNP) drift-diffusion system. Here the energy is given by the Helmholtz free energy, and we calculated the associated $L^2$ gradient flow, constrained by the PNP system. In future work we shall focus on proving well-posedness of the PNP-type model augmented with diffusion and metabolic cost. This task is severely complicated due to the lack of proper regularity estimates and, moreover, due to the lack of minimum principle for the diffusivity tensor.

\end{document}